\documentclass[11pt]{article}
\usepackage[english]{babel}
\usepackage{graphics}
\usepackage{amsfonts}
\usepackage{amsmath}
\usepackage{amssymb}
\usepackage{mathrsfs}
\usepackage{enumerate}
\usepackage{hyperref}
\usepackage{epsfig}
\usepackage{amscd}
\usepackage{amsmath}

\newtheorem{lemma}{Lemma}[section]
\newtheorem{corollary}{Corollary}[section]
\newtheorem{proposition}{Proposition}[section]
\newtheorem{example}{Example}[section]
\newtheorem{theorem}{Theorem}[section]

\newtheorem{remark}{Remark}[section]
\newtheorem{definition}{Definition}[section]
\def\psn{\par\smallskip\noindent}
\def\QED{\hfill $\Box$\par\medskip\noindent}

\def\scatola{\lower5pt\hbox{\vbox{\hrule\hbox{\vrule\kern2pt\vbox%
{\kern5pt\hbox{\mathsurround=0pt }\kern2pt}\kern4pt\vrule}\hrule}}\
} %%
\def\l{\lambda}

\def\t{{\tau}}
\def\noin{\noindent}

\def\H{{\mathbb H}}

\def\R{{\mathbb{R}}}

\def\lev{\mathrm{lev}}

\def\l{\lambda}

\def\0{\overline 0}

\def\psn{\par\smallskip\noindent}

\def\Kca{{\cal K}}
\font\tenmsb=msbm10 \font\sevenmsb=msbm7 \font\fivemsb=msbm5
\newfam\msbfam
\textfont\msbfam=\tenmsb \scriptfont\msbfam=\sevenmsb
\scriptscriptfont\msbfam=\fivemsb

\newcommand*{\finedim}{{\hfill{\lower5pt\hbox{\vbox{\hrule\hbox{\vrule\kern2pt\vbox%
{\kern5pt\hbox{\mathsurround=0pt
}\kern2pt}\kern4pt\vrule}\hrule}}\bigskip }}}

\date{}

\begin{document}

\title{On the $\H$-cone-functions for H-convex sets.}
\author{
 A. Calogero\thanks{
Dipartimento di Matematica e Applicazioni, Universit\`a degli Studi
di Milano-Bicocca, Via Cozzi 55, 20125 Milano, Italy ({\tt
andrea.calogero@unimib.it}, $^+$corresponding author: {\tt
rita.pini@unimib.it})},
 R. Pini$^\dag$$^+$
 }

%\date{}
\maketitle

\begin{abstract}
\noindent
Given a compact and H-convex subset $K$ of the Heisenberg group $\H$, with the origin $e$ in its interior,
we are interested in finding a homogeneous H-convex function $f$ such  that
$f(e)=0$ and $f\bigl|_{\partial K}=1$; we will call this function $f$ the  $\H$-cone-function of vertex $e$ and base $\partial K$.  While the equivalent version of this problem in the Euclidean framework has an easy solution, in our context this investigation turns out to be quite entangled, and the problem can be unsolvable. The approach we follow makes use of an extension of the notion of convex family introduced by Fenchel. We provide the precise, even if awkward, condition required to $K$ so that $\partial K$ is the base of an $\H$-cone-function of vertex $e.$ Via a suitable employment of this condition, we prove two interesting binding constraints on the shape of the set $K,$ together with several examples.
\end{abstract}

%\textbf{ALTRI TITOLI: }
%
%\textbf{Some properties of cone-functions on $\H$}
%



\section{Introduction}
Given a compact and convex set $K$ in $\R^n$, with the origin $e$ in its interior, the problem of constructing the homogeneous  convex function $f: {K}\to [0,+\infty)$ so that $f(e)=0$ and $f\bigl|_{\partial K}=1$ has always a unique solution; this function is called cone-function with vertex $e$ and base $\partial K$, and it is the biggest convex function with these constraints in $e$ and in $\partial K$.
If we move the vertex from the origin to a generic point $x_0$ in the interior of $K,$ and require that the convex function has some positive value on the boundary of $K$, then, by translation and multiplication by a constant, the problem is strictly related to the initial one.

Among the several applications of the cone-functions, our attention is drawn by the Alexandrov Maximum Principle
(see \cite{Gu2001}) and the study of the engulfing property for convex functions in order to obtain some regularity for the solution of the Monge-Amp$\grave{\rm e}$re equation in the spirit of Caffarelli (see \cite{Caf1990}, \cite{Caf1991}). Homogeneity and convexity play both a fundamental role (see Section \ref{section 3}).

In principle, the construction of the cone-function $f$ with vertex $e$ and base $\partial K$ can be undertaken following two equivalent approaches. The first one consists in giving directly the graph of $f$: this is the surface obtained by taking the half lines with the origin  of $\R^{n+1}$ as endpoint, and passing through every point of $\partial K\times\{1\}\subset{\R^{n+1}}.$
Actually, this approach cannot be carried on in the Heisenberg case, due its sub-Riemannian structure.

  The second approach relies on the remark that every function is characterized via its level sets (see the definition below). Thereby, one is lead to consider the family of sets  ${\mathcal F}=\{K_\t\}_{\t\ge 0}$, where $K_\t$ is obtained  by dilating $K$ (i.e., $K_\t=\t K$),  and
 wonder whether the function $f$ defined by the property $f\bigl|_{\partial K_\t}=\t$, for every $\t\ge 0,$ turns out to be convex.

  The study of the families of  sets ${\mathcal F}$ giving rise to a convex function dates back to the pioneering work of Fenchel (see \cite{Fe1953}). Since we will follow this approach in the framework of the Heisenberg group, let us illustrate briefly the main
   ideas in a very general setting.

Given a nonempty set $\Omega,$  a real number $\tau,$ and a function $f:\Omega\to \R,$ we define the level set $\lev_{\le\tau}(f)$ of $f,$ for the level $\tau,$ as follows:
$$
\lev_{\le \tau}(f)=\{x\in \Omega:\ f(x)\le \tau\}.
$$
If $\tau<\inf f,$ then $\lev_{\le \tau}(f)=\emptyset.$ If we denote by $I$ the smallest interval containing the range of $f,$ then the interval $I$ has some interior points, unless $f$ is a constant function.
It is easy to show that the family of subsets $\{\lev_{\le\tau}(f)\}_{\tau\in I}$ enjoys some properties, that we collect in the following
\begin{definition}\label{def level sets family} {A family $\{K_{\tau}\}_{\tau\in I}$ of subsets of $\Omega$ is said to be a {family of level sets} in $\Omega$ if
\begin{itemize}
\item[{\rm I.}] $\displaystyle \bigcup_{\tau\in I}K_\tau=\Omega$;
\item[{\rm II.}] $K_{\tau_1}\subseteq K_{\tau_2},$ if $\tau_1<\tau_2$;
\item[{\rm III.}] $\displaystyle \bigcap_{\tau>\tau_0}K_\tau=K_{\tau_0},$ and $\displaystyle \bigcap_{\tau\in I}K_\tau=\emptyset$ if $I$ is open to the left.
\end{itemize}}
\end{definition}
It is interesting to point out that any family of sets $\{K_{\tau}\}_{\tau\in I}$ of $\Omega$ indexed on an interval $I,$ and satisfying {I.-III.} in the definition above, gives rise to a unique function $f:\Omega\to \R$ via the formula
\begin{equation}\label{eq:function via family}
f(x)=\inf \{\tau\in I:\, x\in
K_\tau\}=\min \{\tau\in I:\, x\in K_\tau\};
\end{equation}
in particular, $\lev_{\le \tau}(f)=K_\tau,$ thereby establishing a one-to-one correspondence between the functions $f$ defined on $\Omega$ and the families of level sets of $\Omega.$

%Let us suppose now that the set $\Omega$ is a nonempty and convex subset of $\R^n.$ If $f:\Omega\to \R$ is a convex function, i.e.,
%$$
%f((1-\theta)x+\theta y)\le (1-\theta)f(x)+\theta f(y),\quad \forall x,y\in \Omega, \theta\in [0,1],
%$$
%then $\lev_{\le \tau}(f)$ is a convex subset of $\Omega,$ for every $\tau.$
As a matter of fact, if $\Omega\subseteq \R^n,$ and ${\mathcal K}=\{K_\tau\}_{\tau\in I}$ is a family of convex subsets of $\Omega,$ then the function $f$ defined via \eqref{eq:function via family} is not convex, in general, but it satisfies the weaker condition
$$
f((1-\theta)x+\theta y)\le \max\{f(x),f(y)\},\quad \forall x,y\in \Omega, \theta\in [0,1],
$$
i.e., $f$ is always a quasi-convex function.
In order to get a convex function starting from a family of level sets,  the assumption of convexity of each set $K_\tau$ is too weak, and the following stronger condition interconnecting the nested sets $\{K_\tau\}$ should be taken into account:
\begin{equation}
(1-\theta)K_{\tau_0}+\theta K_{\tau_1}\subseteq K_{\tau_\theta},\tag{\textbf{C}}\label{eq:condition C}
\end{equation}
for every $\tau_0,\tau_1\in I,$ $\theta\in [0,1],$ $\tau_{\theta}=(1-\theta)\tau_0+\theta \tau_1.$
%Condition \eqref{eq:condition C} trivially implies that every set $L_\tau$ is convex; indeed, by taking $\tau_0=\tau_1=\tau,$ the condition
%$$
%(1-\theta)K_\tau+\theta K_\tau\subseteq K_\tau
%$$
%is equivalent to the convexity of $K_\tau.$
In \cite{Fe1953} Fenchel showed that condition (\textbf{C}) does indeed characterize the sets of a family ${\mathcal F}$ giving rise to a convex function $f$ via
 \eqref{eq:function via family}; such a family was named \emph{convex family in $\Omega$}. He proved the following:
\begin{theorem}\label{theorem condition C}
Let $\Omega$ be a nonempty convex subset of $\R^n$ and let ${\mathcal F}=\{K_\tau\}_{\tau\in I}$ be a family of level sets in $\Omega$. Then, ${\mathcal F}$  is a convex family  if and only if condition \eqref{eq:condition C} is fulfilled.
\end{theorem}
In the particular case ${\mathcal F}=\{\t K\}_{\t\ge 0}$, where $K$ is a compact and convex set with $e$ in the interior,
 the theorem above applies and gives an answer to the problem raised in the beginning of this introduction. As a matter of fact,
the function $f$ defined via \eqref{eq:function via family} is exactly
the cone-function of vertex $e$ and base $\partial K$.

Our aim is to investigate whether, and under what conditions, these ideas can be fruitfully employed in order to construct an $\H$-cone-function, that is, a cone-function in the Heisenberg group $\H$. In this context, we deal with a notion of convexity, the so-called H-convexity, both for sets and functions, which is related to the sub-Riemannian structure of $\H$ (see Definitions \ref{def h convex set}, \ref{def h convex function}). We recall only that  $\R^3$-convexity implies H-convexity.

In this paper, we will follow
the Fenchel approach. First, we study the families  ${\mathcal F}=\{K_\tau\}_{\tau\in I}$ of level sets in $\H$ giving rise to H-convex functions. In the spirit of Theorem \ref{theorem condition C}, we provide an apparently simple but quite unhandy characterization for ${\mathcal F}$ to be an H-convex family: condition  (\textbf{C$_H$}) in Proposition \ref{propo} turns out to be the Heisenberg version of the previous condition (\textbf{C}).
The second step of the paper is to consider a compact, H-convex set $K\subset \H,$ with $e$ in the interior, and therefore deal with the family ${\mathcal F}=\{K_\t\}_{\t\ge 0}$, where $K_\t=\delta_\t K$ and $\delta_\t$ is the non-isotropic dilation of $\H$. Our purpose is to detect what conditions must be fulfilled by $K$ in order that $\mathcal{F}$ is an H-convex family. This function, that will be, by construction,
homogeneous with respect to the dilations $\delta_\l,$ will be called \emph{$\H$-cone-function}, and $K$ will be said to generate an \emph{H-convex family} via dilations.
It is worthwhile noting that in $\H,$ contrary to the Euclidean case, the notion of H-convexity and the dilations $\delta_\l$  have a good mutual behaviour only along the horizontal directions.

Let us devote a few lines to make a delicate point clearer. Suppose that $K$ is compact and $\R^3$-convex, with $e$ in the interior. Then, the cone-function $f$ with vertex $e$ and base $\partial K$ does exist; it is an $\R^3$-convex function, thereby H-convex. Starting from the same set $K$, and working with the non isotropic dilations of $\H,$ it may happen that we can construct the $\H$-cone-function $\phi$ with vertex $e$ and base $\partial K$; in this case, this function $\phi$ is $H$-convex, but obviously it is different from the previous Euclidean cone $f$ (see Examples \ref{example bolla E in R} and \ref{example bolla E in H}). However, it may be that such $\H$-cone-function $\phi$ does not exist (see Examples \ref{example cilinder}).

More generally, let $K$ be H-convex, not $\R^3$-convex; in this case, undoubtedly the cone-function with vertex $e$ and base $\partial K$ does not exist (otherwise, its level sets would be $\R^3$-convex); but it may happen that we can construct an $\H$-cone-function with vertex $e$ and base $\partial K$.

Our interest is essentially motivated by the attempt to approach the problem of the engulfing properties for H-convex functions using the Alexandrov Maximum Principle
(see \cite{CapMal2006}, \cite{BaCaKr2015}). The reader could be surprised that there exist H-convex sets that do not admit $\H$-cone-functions. However, in order to study the engulfing property for H-convex functions, whose level sets are H-convex, this seems to be the right line of investigation.
In this initial study about the existence of $\H$-cone-functions, we try to shed some light on the \lq\lq shape\rq\rq\ of a set $K$ that generates an H-convex family.  Whenever $K_\tau=\delta_\tau K,$ this condition can be written in terms of closedness of $K$ with respect to suitable parabolic arcs with ending points in $K$ (Proposition \ref{prop condizione parabola}).
Taking advantage of this characterization, and with the focus on radial sets, i.e., sets invariant with respect to a rotation around the $t$ axis, we are able to prove our main result, that concern some binding constraints that $K$ must fulfill:
\begin{theorem}\label{ball}
Let $K\subset \H$ be a radial and compact set which generates an H-convex family via dilations. The following statements hold:
\begin{enumerate}[i.]
\item if $\xi_0=(x_0,y_0,0)\in \partial K$, then
$$
B_G(e,\|\xi_0\|_G)\subset K,
$$
where $B_G(\xi,r)$ denotes the Kor\'{a}nyi ball with center $\xi$ and radius $r;$
\item the projection of $K$ onto the $xy$-plane is contained in $K.$
\end{enumerate}
\end{theorem}
%In particular, part ii. of the theorem extends a result recently obtained by Balogh et al. for homogeneous

%\begin{theorem}\label{ball}
%Let $K\subset \H$ be a radial and compact set which generates an H--convex family via dilations. If $\xi_0=(x_0,y_0,0)\in \partial K$, then
%$$
%B_G(e,\|\xi_0\|_G)\subset K,
%$$
%where $B_G(\xi,r)$ denotes the Kor\'{a}nyi ball with center $\xi$ and radius $r.$
%\end{theorem}
%As a consequence, the gauge function turns out to be the biggest radial, H-convex and homogeneous function with prescribed value at any point of the $(x,y)$-plane.
%
%
%\begin{theorem}\label{teorema proiezione}
% Let $K\subset \H$ be a radial and compact set which generates an H-convex family via dilations.
% Then, the projection of $K$ onto the plane  $\{t=0\}$ is contained in the intersection of $K$ with $\{t=0\}$, i.e.,
% $$
%{\rm Pr}_{\{t=0\}}(K)\subset \left(K\cap \{(x,y,0)\in\H\}\right).
%$$
%\end{theorem}

\noindent\textbf{Acknowledgements} We wish to thank Nicolas Hadjisavvas for some initial stimulating suggestions and for drawing our attention on the unpublished work by Fenchel about a family of convex sets.

\section{H-convex sets, families and functions in $\H$}

\subsection{Preliminaries on the Heisenberg group $\H$}

In this short subsection we collect some essential notions concerning the Heisenberg group that will be useful for the sequel. For a thorough discussion on this topic, see, for instance, \cite{CaDaPaTy2007}

The Heisenberg group $\H$ is the Lie group given by the underlying
manifold $\R^3$ with the non commutative group law
$$
\xi\circ \xi'=(x,y,t)\circ(x',y',t')= \left(x+x',y+y',t+t'+2(x'y-xy')\right),
$$
unit element $e=(0,0,0),$ and $\xi^{-1}=(-x,-y,-t)$. Left translations
and anisotropic dilations are, in this setup, $L_{\xi_0}(\xi)=\xi_0\circ
\xi$ and $\delta_\l(x,y,t)=\left(\lambda x,\lambda
y,\lambda^2t\right).$ A function $f:\H\to \R$ is said to be homogeneous if
$$
f(\delta_\l\xi)=\l f(\xi),\qquad \forall \ \xi\in\H,\ \lambda>0.
$$
In the whole paper, the homogeneity, unless otherwise specified, will be related to the dilation $\delta_\l.$

 The  differentiable  structure on $\H$ is
determined  by  the left invariant vector fields
$$
X=\partial_{x}+2y\partial_{t},\qquad
Y=\partial_{y}-2x\partial_{t},\qquad Z=\partial_{t},\qquad
\texttt{\rm with}\ \ [X,Y]=-4Z.
$$
\noin The vector field $Z$ commutes with the vector fields $X$ and
$Y$; $X$ and $Y$ are called \it horizontal vector fields\rm.
The Lie algebra $\mathfrak{h}$ of $\H$ is the stratified algebra
$\mathfrak{h}=\R^3\simeq V_1\oplus V_2,$ where $V_1={\rm
span}\left\{X,Y\right\},$ $V_2={\rm span}\left\{Z\right\};$
$\langle\cdot ,\cdot\rangle$ will denote the inner product in $V_1\simeq\R^2$. Via the
exponential map $\exp:\mathfrak{h}\to\H$ we identify the vector
$\alpha X+\beta Y+\gamma Z$ in $\mathfrak{h}$ with the point
$(\alpha, \beta, \gamma)$ in $\H.$

% the inverse $\nu:
%\H\to\mathfrak{h}$ of the exponential map has the unique
%decomposition $\nu=(\nu_1,\nu_2),$ with $\nu_i:\H\to V_i.$ Since we
%identify $V_1$ with $\R^2$ when needed, $\nu_1:\H\to V_1\cong
%\R^2$ is
%given by $\nu_1(x,y,t)=(x,y):=Pr(x,y,t).$

The sub-Riemannian structure of the Heisenberg group gives rise, among other things, to some \lq\lq horizontal\rq\rq\ notions.
We say that $r$ is a
{horizontal line} through $\xi_0$ if
$r=\{\xi=L_{\xi_0}(\exp(s w))\in\H,\ s\in\R\}$,
for some  $w\in V_1$ fixed.
 We denote by ${\mathcal R}_{\xi_0}$ the set of all the horizontal lines through $\xi_0.$
 %; clearly ${\mathcal R}_{\xi_0}=L_{\xi_0}({\mathcal R}_{e}),$ i.e. every horizontal line $r\in in $.
 A horizontal segment is a convex and bounded subset of a horizontal line.

\noindent Given a point $\xi_0\in\H$, the {horizontal plane}
$H_{\xi_0}$ associated to $\xi_0$ is the plane in $\H$ defined by
$$H_{\xi_0}=L_{\xi_0}\left(\exp(V_1)\right)=\left\{\xi=(x,y,t):
t=t_0+2y_0x-2x_0y\right\}.
$$
Note that  $\xi'\in H_{\xi}$ if and only if $\xi\in H_{\xi'}.$
%Moreover,  $r\in{\mathcal R}_{\xi_0}$ implies $r(s)\in H_{\xi_0}$.

 %For more details on the structure of the Heisenberg group see, e.g., \cite{CaDaPaTy2007} and \cite{BoLaUg2007}.
%
%
%\vskip 0.1truecm The main issue in the analysis on the Heisenberg
%group is that the classical differential operators are considered
%only in terms of the horizontal fields. For any open subset $\Omega$
%of $\H,$ let us denote by $\Gamma^1(\Omega)$ the class of functions
%having continuous derivatives with respect to the vector fields $X$
%and $Y.$ We recall that the horizontal gradient of a function
%$u\in\Gamma^1(\Omega)$ at $\xi\in\Omega$ is the $2$-vector
%$$
%(\nabla_H u)(\xi)= \left((Xu)(\xi), \ (Yu)(\xi)\right),
%$$
%written with respect to the basis $\{X,Y\}$ of $V_1.$

It is well known that a homogeneous norm $N:\H\to[0,+\infty)$ is a homogeneous function satisfying the classical properties of the norm, with respect to the law of the group.
As a particular case, the  Kor\'{a}nyi norm  is defined  by
\begin{equation}\label{gauge}
 \| \xi\|_G=\left((x^2+y^2)^2+t^2\right)^{1/4},\quad \forall \xi=(x,y,t)\in\H.
\end{equation}
Furthermore, we will denote by $\|\cdot\|_q$ the homogeneous quasi-norm defined as follows:
 $\|\xi\|_q=\max\{\|(x,y)\|_E,\, |t|^{1/2}\}$, for every $\xi=(x,y,t)\in\H$.
 %the Euclidean norm in $\R^n$; clear it has a good behaviour w.r.t. the dilation $\l\mapsto \l x,$ for every $x\in\R^n$ and $\l>0$; for this reason we emphasize \lq\lq homogeneous\rq\rq\ when iv. holds.
We will denote by $B_G (\xi_0,r)$  the  closed Kor\'anyi ball with center $\xi_0$ and radius $r\ge0,$ defined by a left translation: more precisely,
$$
B_G (\xi_0,r)=L_{\xi_0}\left(\left\{\xi\in\H: \|\xi\|_G\le r\right\}\right).
$$
Similarly, $B_q(\xi_0,r)$ will denote the closed ball with center $\xi_0$ and radius $r\ge0$ associated to $\|\cdot \|_q.$
We denote by $B_E(\xi_0,r)$ the usual closed Euclidean ball with centre $\xi_0\in\R^n$ and radius $r\ge0,$ and by $\|\cdot \|_E$ the usual Euclidean norm in $\R^n$.
The topology induced in $\H$ by all these balls are equivalent.

\subsection{H-convex sets, families and functions}

Let us start with the definition of H-convex set (see \cite{DaGaNh2003}):

\begin{definition}\label{def h convex set}
A subset $\Omega$ of $\H$ is H-convex if it
contains every horizontal segment with endpoints in $\Omega,$ i.e., if $\xi\in\Omega$ and $\xi'=\xi\circ \exp v\in\Omega$ for some $v\in V_1$, then
$$ \xi\circ\exp (\theta v)\in H_\xi\cap\Omega,\qquad \forall\theta\in[0,1].$$
\end{definition}
It is clear, by the definition, that an $\R^3$-convex set
$\Omega$ is H-convex in the Heisenberg group, since every horizontal
segment is a particular $\R^3$-segment.

The regularity of these H-convex sets was studied in \cite{ArCaMo2011} by means of intrinsic cones at non-characteristic points. Let us devote a few words to
one of the results therein.
For a set $\Omega\subset\H$, a point $\xi_0\in\partial \Omega$ is said to be non-characteristic for $\Omega$ if there exists a horizontal line $r$ through $\xi_0$ so that
$r$ enters $\Omega$ at $\xi_0$, i.e., $r\cap \texttt{\rm int}(\Omega)\cap  B_q(\xi_0,\rho)\not=\emptyset,$ for all $\rho>0$.

For every horizontal line $r\in{\mathcal R}_{e},$ let   $\pi_r:\H\to r$ be the usual (Euclidean) orthogonal projection onto $r,$ and let us define the operator
$\pi_r^\perp:\H\to\H$ to be the $\H$-projection via the identity
$$
\xi=\pi_r^\perp(\xi)\circ \pi_r(\xi),\qquad\forall\xi\in\H.
$$
We emphasize that $\pi_r^\perp$ is strictly related to the non commutative group operation in $\H$.
The open cone $C(e,r,\alpha,h)$ with vertex $e\in\H$, axis  $r\in{\mathcal R}_{e}$, aperture $\alpha>0$, and height $h>0$  is given by
$$
C(e,r,\alpha,h)=\{\xi\in\H:\ \|\pi_r^\perp(\xi)\|_q<\alpha \|\pi_r(\xi)\|_q<\alpha h\}.
$$
The open cone $C(\xi_0,L_{\xi_0}(r),\alpha,h)$ is defined by a left translation, i.e.,
$$C(\xi_0,L_{\xi_0}(r),\alpha,h)=L_{\xi_0}\left(C(e,r,\alpha,h)\right).$$
In \cite{ArCaMo2011}, the authors prove the following result, that provides useful information about the shape of an H-convex set:
\begin{theorem}\label{teo cono Arena Caruso Monti}
Let $\Omega$ be an H-convex set. Let $\xi_0$ be a non-characteristic point for $\Omega$, with $L_{\xi_0}(r)\in {\mathcal R}_{\xi_0}$ entering  $\Omega$ at $\xi_0$. Then, there exist $\rho>0$, $\alpha>0$, and $h>0$, so that for all $\xi\in B_q(\xi_0,\rho)\cap \partial \Omega$ we have
$$
C(\xi,L_{\xi}(r),\alpha,h)\subset \texttt{\rm int}(\Omega).
$$
\end{theorem}

Let us now deal with the notion of H-convexity for functions:
\begin{definition}\label{def h convex function} { Let $\Omega\subset \H$ be H-convex.
A function $f:\Omega\to \R$ is said to be {H-quasiconvex} if and only
if
$$
f(\xi\circ\exp (\theta v))\le \max(f(\xi), f(\xi')),
$$
for every $\xi\in \Omega,$ $\xi'=\xi\circ\exp(v)\in H_\xi\cap \Omega$ and
$\theta\in [0,1];$ equivalently, all the level sets
$\lev_{\le \tau}(f).$

A function $f:\Omega\to \R$ is said to be {H-convex} if
it is $\R$-convex along every horizontal segment within $\Omega$, i.e.,
\begin{equation}\label{def:H-convexity}
f(\xi\circ\exp (\theta v))\le (1-\theta)f(\xi)+\theta f(\xi'),
\end{equation}
 for every $\xi\in \Omega,$ $\xi'=\xi\circ\exp v\in H_\xi\cap \Omega$ and
$\theta\in [0,1].$}
\end{definition}

It follows trivially that the $\R^3$-convexity implies the H-convexity in the Heisenberg group, since every horizontal
segment is a particular $\R^3$-segment; nevertheless, there exist H-convex functions that are not $\R^3$-convex function. An example is provided by the
 Kor\'anyi gauge, i.e., the function $\xi\mapsto \| \xi\|_G$.

The definition of H-convexity implies that the
level sets of an H-convex function are H-convex
sets. However, like in the Euclidean case, if all the levels $\lev_{\le \tau}(f)$ are
H-convex, one cannot infer that the $f$ is H-convex.
Let us consider, for instance, the function $f:\H\to\R$ defined by
$f(x,y,t)= \ln(x^2+y^2+t^2+1),$ and note that all the level
sets are Euclidean convex, in particular H-convex. Since
the horizontal Hessian is not positive semidefinite everywhere, we
deduce that this function is not H-convex (see \cite{CaCaPi2008} for details and a thorough analysis of the H-quasiconvex functions, named \emph{weakly {H}-quasiconvex} function therein).

Likewise the Euclidean framework, given a function $f:\H \to \R$ and a point $\xi_0\in \H,$ one defines the so-called horizontal subdifferential $\partial_H f(\xi_0),$ that is the (possibly empty) subset of $V_1:$
$$
\partial_H f(\xi_0)=\{p\in V_1:\ f(\xi_0\circ\exp(v))\ge f(\xi_0)+\langle p,v\rangle,\ \forall v\in V_1\}.
$$
In \cite{CaPi2011},  it is shown that the nonemptiness of $\partial_H f(\xi)$ at every $\xi\in \H$ completely characterizes   the H-convexity of $f.$ In addition, via an integration technique, the function $f$ can be built up starting from the horizontal subdifferential.

Despite the notion of H-convexity is affected by the horizontal structure, the regularity of H-convex functions is
high. In fact, in \cite{BaRi2003} the authors prove that any H-convex function is locally Lipschitz
continuous with respect to the sub-Riemannian Carnot-Carath\'{e}odory metric. This is weaker than the Euclidean Lipschitz continuity, but is implies that every H-convex function is, at least, continuous.

Motivated by the previous notions of H-convexity and H-quasiconvexity, and following the ideas of Fenchel \cite{Fe1953}, we give the following definition:

\begin{definition}\label{def H-quasiconvex family and H-convex family}
%[H-quasiconvex family and H-convex family]
{Let $\Kca=\{K_\t\}_{\t\in I}$ be a family of level sets in $\H$ (see Definition \ref{def level sets family}).
We say that $\Kca$ is an H-quasiconvex family (H-convex family) in $\H$ if the sets $\{K_\t\}_{\t\in I}$ are the level sets of an H-quasiconvex (H-convex) function $f:\H\to\R$, i.e., $\lev_{\le \t}(f)=K_\t,$ for $\t\in I$.}
\end{definition}
As mentioned in the introduction, there is a unique function $f:\H\to \R$ associated to $\Kca$ as in the previous definition, and it is given by the formula
\begin{equation}\label{eq:function via family in H}
f(\xi)=\inf \{\tau\in I:\, \xi\in
K_\tau\}=\min \{\tau\in I:\, \xi\in K_\tau\}.
\end{equation}
The following remark trivially holds:
\begin{remark}\label{convexity_of_sets}
If $\{K_\t\}_{\t\in I}$ is an H-convex family, then, for every $\t\in I,$ the set $K_\t$ turns out to be H-convex, and closed.
\end{remark}

In the next result condition (\textbf{C$_H$}) will come into play;  this condition forces the nested sets $\{K_{\tau}\}$ to have a mutual \lq\lq convex\rq\rq\ behaviour when a horizontal movement is involved. Despite being quite unhandy, its role will be prominent in the sequel: this is, in fact, the precise condition that an H-convex family must fulfill, just like a convex family in the Euclidean case is identified by Theorem \ref{theorem condition C}.
\begin{proposition}\label{propo} Let $\{K_\t\}_{\t\in I}$ be a family of level sets in $\H,$ and denote by $f$ the related function defined in \eqref{eq:function via family in H}. Then,
$f$ is H-convex
if and only if, for every $\t_0,\t_1\in I,$ for every $v\in V_1,$ and for every $\theta\in [0,1],$
\begin{equation}
K_{\t_0}\cap (K_{\t_1}\circ \exp v)\subseteq K_{\t_{\theta}}\circ \exp (\theta v),\tag{\textbf{C$_H$}}\label{eq:condition CH}
\end{equation}
where $\t_{\theta}=(1-\theta)\t_0+\theta \t_1.$
\end{proposition}
\textbf{Proof:} Let us prove the \lq\lq if\rq\rq\ part. Take any $\xi_0\in \H,$ $v\in V_1,$ and set $\xi_1=\xi_0\circ \exp v.$ Let $f(\xi)=\min_{\xi\in K_\t}\t,$  $\t_0=f(\xi_0),$ and $\t_1=f(\xi_1).$
We have that $\xi_0\in K_{\t_0},$ and $\xi_1\in K_{\t_1}.$ In particular,
$$
\xi_0\in K_{\t_0}\cap (K_{\t_1}\circ \exp (-v)).
$$
From the assumption, $\xi_0\circ\exp (\theta v) \in K_{\t_{\theta}},$ where $\t_{\theta}=(1-\theta)\t_0+\theta \t_1.$ Therefore,
$$
f(\xi_0\circ \exp (\theta v)) =\min_{\xi_0\circ \exp (\theta v)\in K_\t}\t
\le \t_{\theta}
=(1-\theta)f(\xi_0)+\theta f(\xi_1),
$$
thereby proving the inequality \eqref{def:H-convexity}, with $\xi=\xi_0.$

The \lq\lq only if\rq\rq\ part: take $\xi_0\in K_{\t_0}\cap (K_{\t_1}\circ \exp v).$ This implies that $\xi_0\circ\exp(-v)\in K_{\t_1},$ i.e., $f(\xi_0\circ\exp(-v))\le \t_1.$ From
$$
f(\xi_0\circ \exp(-\theta v))\le (1-\theta) f(\xi_0)+\theta f(\xi_0\circ\exp(-v))\le \t_{\theta},
$$
we get that $\xi_0\circ \exp(-\theta v)\in K_{\t_{\theta}},$ and thus the assertion follows.\finedim

\begin{remark}\label{remark worthwhile} Condition \eqref{eq:condition CH} can be restated as follows:
$$
\xi_0\in K_{\t_0},\quad \xi_0\circ \exp v\in K_{\t_1}\Longrightarrow \xi_0\circ \exp (\theta v)\in K_{\t_\theta},
$$
for all $\t_0,\t_1\in I,$ $\theta\in [0,1],$ and $\t_\theta=(1-\theta)\t_0+\theta \t_1.$
\end{remark}

It is easy to prove, but crucial for the sequel, the following:
\begin{remark}\label{remark dilatazioni}
Let $\Kca=\{K_\t\}_{\t\in I}$ be a family of level sets in $\H$.
\begin{enumerate}[a.]
\item
If $\Kca$ satisfies condition \eqref{eq:condition CH}, then $K_\t$ is H-convex and closed, for every $\t\in I;$
\item[b.]
condition \eqref{eq:condition CH} is invariant with respect to the dilations of $\H,$ i.e.,
 $\Kca=\{K_\t\}_{\t\in I}$ satisfies condition \eqref{eq:condition CH} if and only if,
for every $\alpha>0$, the family of level sets $\Kca^\alpha=\{\delta_\alpha K_\t\}_{\t\in I}$ in $\H$ satisfies condition
 \eqref{eq:condition CH}.
 \end{enumerate}
 \end{remark}

\noindent
Note that, setting $\t_0=\t_1$ in  \eqref{eq:condition CH}, we obtain that if
$\xi_0\in K_{\t_0}\cap (K_{\t_0}\circ \exp v)$, i.e.  $\xi_0$ and $\xi_0\circ\exp(-v)$ belong to $ K_{\t_0},$ then $\xi_0\circ\exp (-\theta v)\in K_{\t_0}$ for every $\theta \in[0,1]$. Hence, in order to obtain the H-convexity in $a.$ of Remark \ref{remark dilatazioni}, the fact that $\Kca$ is a family of level sets is not necessary.

\section{H-convex families via dilations and $\H$-cone-functions}\label{section 3}

In this section we focus on the particular families of sets obtained by dilation of a fixed set $K,$ with the purpose of identifying those sets that give rise to a convex (H-convex) family. This problem is strictly related to the existence of the so called cone-functions ($\H$-cone-function). Let us clarify briefly the problem and its motivation. In order to do that, we recall the following  Alexandrov Maximum Principle:

\begin{theorem}\label{teo alexandrov}
Let $\Omega\subset \mathbb R^n$ be bounded, open and
convex.
Let $u\in C(\overline \Omega)$ be convex
 with $u=0$ on $\partial \Omega.$
Then, for every $x_0\in \Omega,$
\begin{equation}\label{elso-Alex}
    |u(x_0)|^n\leq C_n {\rm dist}(x_0,\partial \Omega){\rm
diam}(\Omega)^{n-1} {\mathcal L}^{n}(\partial u(\Omega)),\
\end{equation}
where $C_n>0$ is a constant depending only on the dimension $n.$
\end{theorem}
Clearly, ${\mathcal L}^{n}$ denotes the Lebesgue measure in $\R^n$.

We recall that the subdifferential of the function $u:\R^n\to\R$ at the point $x_0$ is defined as
$$
\partial u(x_0)=\{p\in \R^n:\ u(x)\ge u(x_0)+\langle p,x-x_0\rangle,\ \forall x\in \R^n\},
$$
and the normal mapping of the function $u$ is the multivalued function given by
$$\partial u(A)=\bigcup_{x\in A} \partial u(x),$$
for every set $A\subset\R^n$.
The main idea of the proof of Theorem \ref{teo alexandrov}
can be summarized in these inclusions:
\begin{equation}\label{dim alexandrov}
\partial v(x_0)\subset \partial v(\Omega) \subset \partial u(\Omega).
\end{equation}
The function $v$ is precisely the cone-function with vertex $x_0$ and base $\partial \Omega$. The \lq\lq shape\rq\rq\ of the set $\partial v(x_0)$ contains indeed all the information about ${\rm dist}(x_0,\partial \Omega)$, ${\rm
diam}(\Omega)$ and $|u(x_0)|$ that are needed in \eqref{elso-Alex}.
The first inclusion in \eqref{dim alexandrov} is clear by definition, but it is an equality, in fact, since the cone-function is homogeneous in $\R^n$ with respect to the center $x_0$.
The second inclusion arises from the comparison principle for the normal mapping, which works taking into account that $u(x_0)=v(x_0)$ and $v=u$ on $\partial \Omega$. For all the details, the reader can give a look at the first chapter of \cite{Gu2001}, but we hope that these hints help to clarify the
role of the cone-function, with its enjoyed properties of homogeneity and convexity.

Moving to the Heisenberg case, some of these ideas can be adapted, even if with quite hard proofs, if we replace in \eqref{dim alexandrov}
$\partial u$ with $\partial_H u$. Again, the homogeneity of a possible $\H$-cone-function is crucial, but now the H-convexity of $v$ is now needed, since the second inclusion in \eqref{dim alexandrov} holds only for H-convex functions (see \cite{BaCaKr2015} for all details).

This section is devoted to the construction of these cone-functions and $\H$-cone-functions. Let us consider, first, the simpler Euclidean case, and then pass to the Heinseberg setting.
In the following, we will assume that $K\subset \H=\R^3$ satisfies the condition
\begin{itemize}
\item[\bf{a.}] $K$ is compact and $e\in \mathrm{int}(K).$
\end{itemize}

\subsection{The Euclidean case}

The problem of constructing a cone-function for a compact and convex set $K$ has a simple solution; the following result holds:

\begin{proposition}\label{Kt dilatazioni euclidee}
%[Euclidean convex families via Euclidean dilations]
Suppose that $K\subset\R^3$ is a convex set satisfying condition \textbf{a.}
Then,  the family $\Kca=\{K_\t\}_{\t\in [0,+\infty)}$, defined by
$$
K_\t=\t K,
$$
is an $\R^3$-convex family.
Moreover, the (convex) function $f:\R^3\to [0,+\infty)$ associated to $\Kca$ by \eqref{eq:function via family} turns out to be homogeneous w.r.t. the Euclidean dilations, i.e., $f(\l \xi)=\l f(\xi),$
 for every $\l\ge 0$ and $\xi\in\R^3$.
\end{proposition}
\textbf{Proof:} First, let us prove that $\Kca$  is a family of level sets in $\R^3$. It is trivial to show that \textbf{a.} implies I. The convexity of $K$ implies II. The proof of III. will be presented hereafter, when dealing with the Heisenberg setting.

   Lets us now show that \eqref{eq:condition C} is fulfilled. Take $\tau_0,\tau_1\in [0,+\infty),$ and $\theta\in [0,1],$ and set $\tau_\theta=(1-\theta)\tau_0+\theta\tau_1.$ Let $\xi\in (1-\theta)\tau_0 K+\theta \tau_1 K;$ then, there exist $\xi_1,\xi_2\in K$ so that
$$
\xi=(1-\theta)\tau_0 \xi_1+\theta \tau_1 \xi_2=\tau_{\theta}\left(\frac{(1-\theta)\tau_0}{\tau_\theta}\xi_1+\frac{\theta\tau_1}{\tau_\theta}\xi_2\right).
$$
Since $\xi=\tau_{\theta}\xi',$
where $\xi'$ is a convex combination of points of $K,$ we conclude that $\xi\in \tau_\theta K.$ It is an easy exercise to prove that $f$ is homogeneous with respect to the Euclidean dilations.

\QED

We will say that the function $f$ in Proposition \ref{Kt dilatazioni euclidee} is the cone-function
with vertex $e$ and base $\partial K$. This definition is motivated by the shape of the graph of $f$ in $\R^4$, that is a cone
with vertex $(e,0)\in\R^3\times\R$ and base $\partial K,$ having \lq\lq value\rq\rq\ 1 on $\partial K$.

The following simple example of cone-function in $\R^3$ is interesting when compared to the forthcoming Example \ref{example bolla E in H} in $\H,$ where we start from the same set $K,$ but we deal with different dilations:
\begin{example}[$K$ unit ball and Euclidean dilations]\label{example bolla E in R}\rm{
Let $K=B_E (0,1).$ Trivially, $K$ satisfies the assumptions of Proposition \ref{Kt dilatazioni euclidee}; $K_\t=\t K=B_E(0,\t)$ for every $\t>0,$ and the cone-function $f:\R^3\to\R$  defined by \eqref{eq:function via family} is $f(x,y,t)=\|(x,y,t)\|_E$.}
\end{example}

\subsection{The Heisenberg  case}

In this subsection we deal with the problem of constructing an H-convex family arising by dilations of a fixed set $K\subset \H,$ and, as a byresult, identifying a cone-function in $\H$ with vertex $e$ and base $\partial K.$
Let us now consider the  family of sets $\Kca=\{K_\t\}_{\t\in [0,+\infty)}$ defined by
\begin{equation}\label{delta K}
K_\t=\delta_\t K,\qquad \t\ge 0;
\end{equation}
For these type of families $\Kca$ the following proposition holds:
\begin{proposition}
Suppose that $K$ satisfies assumption \textbf{a.} and
\begin{itemize}
\item[\bf{b.}] $\delta_\t K\subset K,$ for every $\t\in [0,1).$
\end{itemize}
Then, the family $\Kca=\{K_\t\}_{\t\in [0,+\infty)},$  defined via dilations of the set $K$  as in \eqref{delta K}, is a family of level sets in $\H$.
Moreover, if $f:\H\to [0,+\infty)$ is the function associated to $\Kca$ by \eqref{eq:function via family in H}, then
$f$ is homogeneous.
\end{proposition}
\textbf{Proof:} First, let us prove that $\Kca$  is a family of level subsets of $\H$. Again  \textbf{a.} implies I.  Now, take $\t_0<\t_1,$ and $\xi\in \delta_{\t_0}K,$
i.e., $\xi=\delta_{\t_0}\xi'$ for some $\xi'\in K.$ Since $\delta_{\t_0}\xi'=\delta_{\t_1}(\delta_{\t_0/\t_1}\xi'),$ it follows from  \textbf{b.} that $\delta_{\t_0/\t_1}\xi'\in K,$ therefore II. holds.
In order to prove III., note that, from II., $\cap_{\t>\t_0}\delta_{\t}K\supseteq \delta_{\t_0}K.$ To prove the opposite inclusion, let us now take any $\xi$ in the closed set $\cap_{\t>\t_0}\delta_\t K.$
Since $\xi=\delta_\t(\delta_{1/\t}\xi),$ this implies that $\delta_{1/\t}\xi\in K,$ for every $\t>\t_0.$ By the closedness of $K$, and by the continuity of the map $\t\mapsto \delta_{1/\t}\xi,$ we have that $\lim_{\t\to \t_0}\delta_{1/\t}\xi=\delta_{1/\t_0}\xi\in K,$ therefore, $\xi\in \delta_{\t_0}K.$

Denote by $f:\H\to \R$ the function so that $\mathrm{lev}_{\le \tau}(f)=\delta_\tau K.$
Then,
$$
f(\delta_\l\xi)
=\inf\{\tau\in [0,+\infty):\, \xi\in K_{\tau/\l}\}
=\l \inf\{\tau \in[0,+\infty):\, \xi\in K_\tau\}
=\l f(\xi),
$$
thereby proving that $f$ is homogeneous.
\QED

The proposition above and its proof highlight that, while in the Euclidean case the assumption of convexity of $K$
gives II. in Definition \ref{def level sets family}, in the Heisenberg setting the requirement of H-convexity has nothing to do with the dilations, and hence with II. Moreover, as we will see, the H-convexity assumption on $K$ does not guarantee that $\Kca=\{\delta_\tau K\}_{\tau\in [0,+\infty)}$ is an H-convex family (see Example \ref{example cilinder}).

Let us point out a property that is fulfilled by any homogeneous function on $\H$:
 \begin{proposition}\label{subgradient homogenee} Let $f:\H\to [0,+\infty)$ be a  homogeneous function.
If $p\in \partial_Hf(\xi),$ then $p\in \partial_Hf(\delta_\l\xi),$ for every $\lambda>0.$ Moreover, for every set $K$ such that $e\in K$ we have
$$
\partial_H f(K)=\partial_H f(e).
$$
\end{proposition}
\textbf{Proof:} Suppose that $p\in \partial_Hf(\xi).$ Then,
$$
f(\delta_\l\xi\circ \exp v)=f(\delta_\l(\xi\circ \exp(v/\lambda))
=\l f(\xi\circ \exp(v/\lambda))
\ge \l (f(\xi)+\langle p,v/\l\rangle)
=f(\delta_\l \xi)+\langle p,v\rangle.
$$
The second part of the claim follows easily from the upper semicontinuity of the multivalued function $\partial_H f:\H \rightrightarrows V_1$ (see \cite{CaPi2011} for details).
\QED

With the aim of generating H-convex families via dilations of a set $K\subset\H$, we will assume that the set $K$ satisfies conditions \textbf{a.}, \textbf{b.}, and
\begin{itemize}
\item[\textbf{c.}] $K$  is an H-convex set of $\H.$
\end{itemize}

An interesting and natural question is the following:
suppose that $K$ satisfies properties \textbf{a.}-\textbf{c}.; is there any chance that the family $\Kca$ defined  by \eqref{delta K} is an H-convex family?
In this line of investigation we give the following definition:
\begin{definition}\label{K generates an H-convex family via dilations}
%[$K$ generates an H-convex family via dilations]
We say that a compact set $K\subset\H$ generates an H-convex family via dilations, if the family $\Kca$ defined by \eqref{delta K} is an H-convex family. In this case, we say that the function $f:\H\to [0,+\infty)$ associated to $\Kca$ by \eqref{eq:function via family in H} is the  $\H$-cone-function
with vertex $e$ and base $\partial K$.
\end{definition}
Let us emphasize that the $\H$-cone-functions defined above are homogeneous w.r.t. the dilations $\delta_\l,$ and are different from the cone-functions of the Euclidean case (compare, for instance, Example \ref{example bolla E in R} and Example \ref{example bolla E in H}). Moreover, let us mention that in \cite{BiBaMa2009} the authors introduce, for other purposes with respect to our investigations, two different notions of cone-function.

\begin{remark}\label{dopo def}
If a compact set $K\subset\H$ generates an H-convex family via dilations, then $K$ is H-convex and $e\in \mathrm{int}(K).$
\end{remark}
\textbf{Proof:} The set $K$ is H-convex, since $K=\lev_{\le 1}(f)$ and $f$ is an H-convex function. Suppose, by contradiction, that $e\not\in \mathrm{int}(K);$ then, there exists a sequence $\{\xi_n\}_{n}\subset\H\setminus K$, with $\xi_n\to e$. If $\xi_n\not\in K,$ then
$$\lim_{n\to\infty} f(\xi_n)\ge 1;$$ however, since $f$ is H-convex (by the result in \cite{BaRi2003}), it is continuous, and $f(e)=0.$ This leads to a contradiction.
\QED

Let us now spend a few words on the compactness assumption in Definition \ref{K generates an H-convex family via dilations}:  if a family $\Kca$ defined by $K$ via \eqref{delta K} is an H-convex family, then the compactness of $K$ is not an outcome. For example, the set $K=\{(x,y,t)\in\H:\ |x|\le 1\}$ generates an H-convex family via dilations, and the unique function associated to such family via \eqref{eq:function via family in H} is the $\R^3$-convex, H-convex and homogeneous function $f(x,y,t)=|x|$.
However, since we are interested in the construction of cone-functions having as a base the boundary of a compact set, we require in Definition \ref{K generates an H-convex family via dilations} that $K$ is compact.

Let us provide two examples of $\H$-cone-functions:
\begin{example}[$K$ unit  Kor\'anyi ball]\label{example bolla G}\rm{
Let $K=B_G (0,1);$ trivially, $K$ generates an H-convex family, since the function $f(\xi)=\|\xi\|_G,$ associated to $\{\delta\tau K\}$ is H-convex and homogeneous, thereby it is the $\H$-cone-function of $K.$}
\end{example}

\begin{example}[$K$ unit Euclidean ball]\label{example bolla E in H}
\rm{Let $K=B_E (0,1)$.
For every $\tau>0,$ we have
\begin{equation}\label{example bolla 0}
K_\t=\delta_\t K=\left\{(x,y,t)\in\H:\ \frac{x^2}{\t^2}+\frac{y^2}{\t^2}+\frac{t^2}{\t^4}\le 1\right\}
\end{equation}
Denote by $f$ the function defined via \eqref{eq:function via family in H}; clearly, $f(\xi)=\tau$ if and only if $\xi\in\partial K_\tau.$ Since $\xi=(x,y,t)\in \partial K_\tau$ if and only if
$
\t^4-\t^2({x^2}+{y^2})-1=0,
$
an easy calculation gives
\begin{equation}\label{example bolla}
f(x,y,t)=\frac{1}{\sqrt{2}}\sqrt{x^2+y^2+\sqrt{(x^2+y^2)^2+4t^2}}.
\end{equation}
We will show that $f$ is H-convex. In order to do this, we recall first a property of the composition that preserves convexity. Let us consider a function $F:\R^2_+\to [0,+\infty)$, $\R^2$-convex and monotone in the sense that
$$
x_1\le x_1'\ \texttt{\rm and}\ x_2\le x_2'\qquad\Longrightarrow\qquad F(x_1,x_2)\le F(x_1',x_2'),
$$
and the functions $f_i:\H\to [0,+\infty),\ i=1,2,$ both H-convex.
For every $\xi_0\in\H,\ v\in V_1$ and $s\in [0,1],$ we have
\begin{eqnarray*}
&&F\Bigl(f_1(\xi_0\circ\exp(sv)),\,f_2(\xi_0\circ\exp(sv))\Bigr)\\
&&\qquad \le \ F\Bigl((1-s)f_1(\xi_0)+sf_1(\xi_0\circ\exp(v)),\, (1-s)f_2(\xi_0)+sf_2(\xi_0\circ\exp(v))\Bigr)\\
&&\qquad = \
F\Bigl((1-s)(f_1(\xi_0),\ f_2(\xi_0))\,+\,  s(f_1(\xi_0\circ\exp(v)),\, f_2(\xi_0\circ\exp(v))\Bigr)\\
&&\qquad \le \
(1-s)F\Bigl(f_1(\xi_0),\ f_2(\xi_0)\Bigr)\,+\,  sF\Bigl(f_1(\xi_0\circ\exp(v)),\, f_2(\xi_0\circ\exp(v))\Bigr)
\end{eqnarray*}
Since the function $f$ in \eqref{example bolla} is so that
$$
f(\xi)=\frac{1}{\sqrt 2}
F(f_1(\xi),\, f_2(\xi_2)),
$$
where $F(x_1,x_2)=\|(x_1,x_2)\|_E, \quad f_1(x,y,t)=\|(x,y)\|_E\quad \texttt{\rm and}\quad f_2(x,y,t)=\|(x,y,2t)\|_G,$ and $F$ is monotone on $\R^2_+,$ $f_1$ is $\R^3$-convex and, thereby, H-convex, and $f_2$ is H-convex (see Example \ref{example bolla G}), then, by the previous argument, we get that  $f$ is H-convex. Therefore, $K=B_E (0,1)$ generates an H-convex family via dilations.
In particular, $f$ is the $\H$-cone-function with vertex $e$ and base the boundary of the Euclidean ball $B_E (0,1)$.
Finally, we mention that in Example 5.2 in \cite{BalFasSob2018} the authors prove that $f$ is a homogeneous norm, called Lee-Naor norm.
  }
\end{example}

\bigskip Let us go back to the initial question about the construction of the $\H$-cone-function. In order to provide an answer, we write more explicitly condition \eqref{eq:condition CH} in the particular case of families obtained by dilations:

\begin{proposition}\label{prop condizione parabola}
Let $K\subset\H$ such that assumptions \textbf{a.}-\textbf{c} hold, and let us consider the family $\Kca=\{K_\t\}_{\t\in[0,+\infty)}$
defined by $K$ via dilations as in \eqref{delta K}.
Then, $\Kca$ satisfies condition \eqref{eq:condition CH} if and only if,
for every $\xi_0,\ \t>0$ and $(\alpha, \beta)$ so that
\begin{equation}\label{parabola 1 zero}
\xi_0=(x_0,y_0,t_0)\in K,\qquad
\left(\frac{x_0+\alpha}{\t},\frac{y_0+\beta}{\t},\frac{t_0-2(x_0\beta-y_0\alpha)}{\t^2}\right)\in K,\end{equation}
the curve  $s\mapsto\tilde \xi(\theta(s))$ takes its values within $K$, where $\tilde \xi(\theta(\cdot))$ is defined as follows:
 if we set $A=x_0\beta-y_0\alpha,$
\begin{enumerate}[i.]
\item
if $x_0(\t-1)\neq \alpha,$ then  $s\in [{x_0}, \frac{x_0+\alpha}{\t}],$ and
%\begin{equation}\label{parabola 2i}
$$
\tilde \xi(\theta(s))= \left(
s,\,\frac{(y_0(\t-1)-\beta)s+A}{(\t-1)x_0-\alpha},\, \frac{t_0((\t-1)s-\alpha)^2-2A(x_0-s)((\t-1)s-\alpha)}{((\t-1)x_0-\alpha)^2}\right);
%\end{equation}
$$
\item if $x_0(\t-1)=\alpha$ and $y_0(\t-1)\neq \beta,$ the,  $s\in [{y_0}, \frac{y_0+\beta}{\t}],$ and
$$
\tilde \xi(\theta(s))=\left(
{x_0},\, s,\,\frac{t_0((\t-1)s-\beta)^2-2x_0(\beta-(\t-1)y_0)(y_0-s)((\t-1)s-\beta)}{((\t-1)y_0-\beta)^2}\right);
$$
\item if $x_0(\t-1)=\alpha$ and $y_0(\t-1)=\beta,$ then, in this case, $A=0,$ and we get, for $s\in [0,1],$
$$
\tilde \xi(\theta(s))=\left({x_0},\,{y_0},\,\frac{t_0}{(1+s(\t-1))^2}\right).
$$
\end{enumerate}
\end{proposition}
\textbf{Proof:} First of all, due to Remark \ref{remark dilatazioni}-b.,  we can set $\t_0=1$ in \eqref{eq:condition CH}.
%
% In fact
%let us suppose that
%\begin{equation}\label{dim parabola 1}
%\forall \xi_0,\ v,\ \t>0:\quad \xi_0\in K\cap(K_{\t}\circ \exp v)\qquad \Longrightarrow\qquad \xi_0\in K_{(1-\theta+\t\theta)}\circ\exp (\theta v);
%\end{equation}
%let us prove that  \eqref{eq:condition CH} holds for generic $\t_0$ and $\t_1$.  Hence let us suppose that
%$\xi\in K_{\t_0}\cap(K_{\t_1}\circ \exp w);
%$
%setting $\xi_0=\delta_{\frac{1}{\t_0}}(\xi)$, the previous inclusion is, since $K_s=\delta_s(K)$,
%$$
%\xi_0\in \delta_{\frac{1}{\t_0}}\left(K_{\t_0}\cap(K_{\t_1}\circ \exp w)\right)=
%K\cap\left(K_{\t_1/\t_0}\circ \exp (w/\t_0))\right).
%$$
%Now \eqref{dim parabola 1} gives
%$$
%\xi\in
%K_{(1-\theta+\t_1\theta/\t_0)}\circ\exp (\theta w/\t_0)  = \delta_{\t_0}\left(K_{(\t_0(1-\theta)+\t_1\theta)}\circ\exp (\theta w)  \right),$$
%and hence $\xi\in K_{\t_\theta }\circ\exp (\theta w).$
Let us now consider
$
\xi_0\in K\cap\left(K_\t\circ \exp v\right),
$
where $v=(-\alpha, -\beta).$ This is equivalent to
$\xi_0\in K$ and $\delta_{1/\t}\left(\xi_0\circ\exp(-v)\right)\in K,$ that is exactly \eqref{parabola 1 zero}. As in Remark
\ref{remark worthwhile}, condition  \eqref{eq:condition CH} can restated as:
%\begin{equation}\label{parabola 1}
$$
\tilde \xi(\theta)=\delta_{1/\tau_\theta}(\xi_0\circ \exp(\theta v))=\left(\frac{a+\theta \alpha}{\tau_\theta},\frac{b+\theta \beta}{\t_\theta},
\frac{c-2(a\beta-b\alpha)\theta}{\tau_\theta^2}\right)\in K,\quad \forall \theta\in [0,1],
%\end{equation}
$$
where $\t_\theta=1-\theta(1-\t).$
By a change of parametrization, setting $A=a\beta-b\alpha,$ the curve above can be expressed, after tedious computations, as in i.-iii.

\QED

The following example shows that
the problem of finding the $\H$-cone-function can be unsolvable. More precisely, it shows that there exist sets $K$ satisfying \textbf{a.-c.} that do not fulfil \eqref{eq:condition CH}; thus, there are no $\H$-cone-functions with vertex $e$ and the boundary of these sets as bases.

\begin{example}[$K$ cilinder]\label{example cilinder}
\rm{Let $K=\{(x,y,t)\in\H:\ \|(x,y)\|_E\le 1,\ |t|\le 1\}.$
Clearly, $K$ satisfies \textbf{a.-c.}, and it is $\R^3$-convex; nevertheless, \eqref{eq:condition CH} does not hold.
In fact, let us fix $(x_0,\beta)\in\R^2,$ with  $x_0>0,$ $\beta<0 $  and $\| (x_0,\beta)\|_E< 1,$ and consider
$$
\xi_0=(x_0,0,1),\qquad -v=(0,\beta),\qquad \t=\sqrt{1-2x_0\beta}>1.
$$
Under these choices, assumption \eqref{parabola 1 zero} becomes
$$
\xi_0=
(x_0,0,1)
\in K,\quad \mathrm{and}\quad
\left(
\frac{x_0}{\t},\frac{\beta}{\t},1\right)\in K,
$$
and it is obviously fulfilled.
Since $x_0(\t-1)\neq 0$, by Proposition \ref{prop condizione parabola}-i., condition \eqref{eq:condition CH} requires that
$$
\tilde \xi(\theta(s))= \left(
s,\,\frac{\beta(x_0-s)}{(\t-1)x_0},\, \frac{s^2(\t-1+2x_0\beta)-2x_0^2\beta s}{(\t-1)x_0^2}\right)\in K;
$$
for $s\in [x_0/\t, x_0].$  Taking into account that $2x_0\beta=1-\t^2,$ we obtain
$$
\frac{(\t-1+2x_0\beta)}{(\t-1)x_0^2}=-\frac{\t}{x_0^2}<0,
$$
hence the third component of the curve $\tilde \xi(\theta(s))$ is greater than 1 for every $s\in(x_0/\t, x_0).$ This implies that $\tilde \xi(\theta(s))\not\in K$.}
\end{example}

\section{$\H$-cone-functions with base on the boundary of radial H-convex sets}

In order to better understand the role played by condition \eqref{eq:condition CH} in relation to the shape of $K$, and with the aim of constructing an $\H$-cone-function for a compact, H-convex set $K$ having $e$ in the interior, let us focus on a particular class of sets, the radial ones. We say that
$K\subset \H$ is radial if
$$
\xi_0=(x_0,y_0,t_0)\in K\quad \Longleftrightarrow\quad (x,y,t_0)\in K,\ \forall\, (x,y)\in \R^2:\  \|(x,y)\|_E=\|(x_0,y_0)\|_E.
$$
The section investigates, in several steps, the radial sets that are candidates to have an $\H$-cone-function associated.

It is clear that, if $K$  is a compact and radial set, then the $\H$-cone-function $f$ with vertex $e$ and base $\partial K$ (when it exists) is radial, by default, i.e., $f(x,y,t)=g(\sqrt{x^2+y^2},t)$ for some function $g:[0,+\infty)\times \R\to\R$ (equivalently,
 $f(x,y,t)=f(O(x,y),t)$ for every $O\in O(2)$).
The H-subdifferential $\partial_H f$ of these kind of functions enjoys nice properties, as we show:

\begin{proposition}\label{pr:rotation} Let $f:\H\to \R$ be a radial function.
If $p\in \partial_Hf(x,y,t),$ then $Op\in \partial_Hf(O(x,y),t).$

\end{proposition}
\textbf{Proof:} Denote by $\pi:\H\to \R^2$ the projection $\pi(\xi)=z,$ if $\xi=(z,t),$ and by $J:\R^2\to \R^2$ the operator defined as $J(x,y)=(y,-x).$ For every $p\in \partial_Hf(\xi),$ and $O\in O(2),$ we have, for every $v\in V_1\simeq \R^2,$
\begin{align*}
f(Oz+v,t-2\langle Oz,Jv\rangle)&=f(O(z+O^Tv),t-2\langle Oz,JOO^Tv\rangle)\\
&=f(O(z+O^Tv),t-2\langle Oz,OJO^Tv\rangle)\\
&=f(O(z+O^Tv),t-2\langle z,JO^Tv\rangle)\\
&=f(z+O^Tv,t-2\langle z,JO^Tv\rangle)\\
&\ge f(z,t)+\langle z,O^Tv\rangle\\
&=f(Oz,t)+\langle Oz,v\rangle,
\end{align*}
thereby proving the claim.\QED

By collecting the results of the previous proposition and Proposition \ref{subgradient homogenee}, we have the following
\begin{corollary} Let $f:\H\to [0,+\infty)$ be the $\H$-cone-function with vertex $e$ and base $\partial K$, where $K\subset \H$ is radial.
If $p\in \partial_Hf(\xi),$ then $Op\in \partial_Hf(O(\delta_\l\xi)),$ for every $\lambda>0$ and $O\in O(2).$
\end{corollary}

Let us start our investigation on the shape of the radial set $K$ with the following
\begin{proposition}\label{punti non caratteristici}
Let $K \subset \H$ be a radial and compact set which generates an H-convex family.
Let us define $t_{\min}$ and $t_{\max}$ as follows:
$$
t_{\min}=\min\{t:\ (0,0,t)\in K\},\qquad t_{\max}=\max\{t:\ (0,0,t)\in K\}.
$$
If $\xi_0=(x_0,y_0,t_0)\in\partial K$ and $t_0\in(t_{\min},t_{\max})$, then $\xi_0$ is a non-characteristic point.
\end{proposition}
\textbf{Proof:} From Remark \ref{dopo def}, the set $K$ is H-convex, and $e\in \mathrm{int}(K),$ i.e., there exists $\rho>0$ so that $B_G(e,\rho)\subset K.$ In particular, $t_{\min}<0<t_{\max}$.

First, let us suppose that $\xi_0=(x_0,0,0)\in \partial K$, with $x_0>0$. It is clear that the horizontal segment with endpoints $e$ and $\xi_0$ is contained in $K,$ but this does not imply that $\xi_0$ is non-characteristic.
Let $t$ small enough, so that $\xi_t=(0,-t/2,x_0t)\in B_G(e,\rho).$ Since $\xi_t=\xi_0\circ\exp(-x_0,-t/2)\in H_{\xi_0}$, then the H-convexity of $K$ implies that, for every $s\in [0,1],$
$$
\xi_0\circ\exp(s(-x_0,-t/2))=((1-s)x_0,-st/2, x_0ts)\in K.
$$
Furthermore, $K$ is radial; then, for every $s\in[0,1]$ and $\theta\in[0,2\pi),$
 $$
 (r_s\cos\theta,r_s\sin\theta, x_0ts)\in K,\qquad r_s=\|((1-s)x_0,-st/2)\|_E;
 $$
in particular, $(r_s,0,x_0ts)\in K$. Note that $(0,0,x_0ts)\in K,$ because $\|(0,0,x_0ts)\|_G<\rho$; hence, all the horizontal segments with endpoints  $(r_s,0,x_0ts)$ and $(0,0,x_0ts)$ belong to $K$, for every $s\in[0,1]$. The radiality of $K$ entails that, if we consider the horizontal half line $r$ through $\xi_0$
\begin{equation}\label{retta r}
r=\{((1-s)x_0,0,0)\in\H:\ s\ge 0\},
\end{equation}
then it enters $K$, thereby showing that $\xi_0$ is
 non-characteristic.

Now let us consider the case $\xi_0=(x_0,0,t_0)\in \partial K$, with $x_0>0$ and $t_0\in(t_{\min},t_{\max})$. Since $K$ is radial and H-convex, then  $(0,0,t_0)\in K$. Hence, there exists $\rho>0$ so that $B_G((0,0,t_0),\rho)\subset K$. A similar argument as before applies, and proves the assertion.

\QED

\noindent
The previous result reveals that a radial set generating an H-convex family enjoys very nice regularity properties, and we are in the position to apply Theorem \ref{teo cono Arena Caruso Monti} to the points $\xi_0\in\partial K$, where $K$ and $\xi_0$ are as in Proposition \ref{punti non caratteristici}. In  what follows, we will see that a sharper regularity result can be proved for such $K$. To do this, we need the next
\begin{lemma}\label{lemma solido di rivoluzione}
Let $K\subset \H$ be a radial set so that the family $\Kca=\{K_\t\}_{\t\in[0,+\infty)}$, defined by $K$ via dilations as in \eqref{delta K}, satisfies  condition \eqref{eq:condition CH}. If $\xi_0=(x_0,y_0,t_0)\in K$ with $t_0>0$ and $r_0=\|(x_0,y_0)\|_E>0$,
then $K$ contains the solid of revolution $S_{\xi_0},$ defined as follows:
$$
S_{\xi_0}=\left\{(x,y ,\varphi(r))\in\H:\, \ r\in [r_0/\sqrt{2},r_0],\ \|(x,y)\|_E\le r\right\},
$$
where $\displaystyle \varphi(r)=\frac{1}{2}\left(t_0+\sqrt{t_0^2+4r_0^2r^2-4r^4}\right)$.
\end{lemma}
\textbf{Proof:} Without loss of generality, let $\xi_0=(x_0,0,t_0)\in K,$ with $x_0$ and $t_0$ positive, and fix $\theta\in (0,\pi/2)$.
The point $(x_0\cos \theta,x_0\sin \theta ,t_0)$ can be obtained as $(\delta_{\tau}\xi_0)\circ\exp v$ by taking, in
 \eqref{parabola 1 zero},
\begin{equation}\label{espressione tau}
\tau=\frac{\sqrt{x_0^4\sin^2\theta+t_0^2}-x_0^2\sin\theta}{t_0},\qquad  -v=(\alpha,\beta)=(x_0(\tau\cos\theta -1),x_0\tau\sin\theta).
\end{equation}
In this case, $\t<1$ and the curve $\widetilde\xi$ can be expressed as in Proposition \ref{prop condizione parabola}-i.:
\begin{eqnarray*}
&&\widetilde\xi (s)= \Biggl(s,\frac{(x_0-s)\sin\theta}{(1-\cos\theta)},\\
&&\qquad\qquad\quad \frac{t_0((\t-1)s-x_0(\t\cos\theta-1))^2-2x_0^2\t\sin\theta(x_0-s)((\t-1)s-x_0(\t\cos\theta-1))}{\t^2x_0^2(1-\cos\theta)^2}
\Biggr),
\end{eqnarray*}
for $s\in [x_0\cos\theta,x_0].$
Setting $\overline s=x_0(\cos\theta+1)/2,$  the vertex of the parabola $\widetilde\xi$ is given by
$$
(x_V,y_V,t_V)=
\left(\frac{x_0(\cos\theta+1)}{2},\frac{x_0\sin\theta}{2}, \frac{(\tau +1)(t_0(\tau+1)-2x_0^2\tau\sin\theta)}{4\tau^2}\right).
$$
Let us concentrate on $t_V$: a tedious calculation, that takes into account the expression of $\tau$ in \eqref{espressione tau}, gives
$$
t_V=\frac{1}{2}\left(t_0+\sqrt{t_0^2+x_0^4\sin^2\theta}\right).$$
Note that $t_V>c$.
Condition \eqref{eq:condition CH} implies that $(x_V,y_V,t_V)\in K$ for every $\theta$. Furthermore, since $K$ is radial and  $\|(x_V,y_V)\|_E=x_0\sqrt{\frac{1+\cos\theta}{2}}$, we have that
$$
\left(x_0\sqrt{\frac{1+\cos\theta}{2}}\cos \theta',x_0\sqrt{\frac{1+\cos\theta}{2}}\sin\theta', \frac{t_0+\sqrt{t_0^2+x_0^4\sin^2\theta}}{2}\right)\in K,
$$
for every $\theta',$ and $\theta\in(0,\pi/2)$. The previous set of points can be described as a surface of revolution, specifically, as surface obtained by the rotation of a curve $r\to (r,0,\varphi(r))$ around the $t$-axis, where $\varphi:[x_0/\sqrt{2},x_0]\to\R.$ With a change of variable, we obtain
$$
\varphi(r)=\frac{1}{2}\left(t_0+\sqrt{t_0^2+4x_0^2r^2-4r^4}\right).
$$
Hence,
$$
\left\{(x,y ,\varphi(r))\in\H:\,\ r\in [r_0/\sqrt{2},r_0],\ \|(x,y)\|_E= r\right\}\subset K;
$$
since $K$ is H-convex (see Remark \ref{remark dilatazioni}-a.), we have the claim.

\QED

The previous lemma provides an \lq\lq inner constraint\rq\rq\ at a point of the boundary of a set $K$ that generates an H-convex family via dilation. One may wonder whether it is possible to provide global information on this set $K.$ Let us start with an example:

\begin{example}[$K$ cilinder with a Kor\'anyi ball as hat]\label{cilinder with hat}
\rm{Let
$$K=\{(x,y,t)\in\H:\ \|(x,y)\|_E\le 1,\ t^2\le 15\}\cap B_G(e,2).$$
Let us verify that the set above generates an H-convex family via dilations. Let $\xi=(x,y,t),$ and consider the function $f:\H\to\R$ defined as follows:
$$
f(\xi)=\max\left(\|(x,y)\|_E,\ \frac{1}{2}\|(x,y,t)\|_G\right).$$
This function is H-convex, as the maximum of two H-convex functions. Moreover, it is homogeneous, and $\lev_{\le 1}(f)=K$. Hence, $K$ generates an H-convex family via dilation, and $f$ is the $\H$-cone-function of vertex $e$ and base $\partial K$.}

\end{example}

It is clear that the $\H$-cone-function $f$ in the Example \ref{cilinder with hat} is a homogeneous norm. In \cite{BalFasSob2018} the authors prove that, if $N$ is a homogeneous norm, then $N(x,y,t)\le N(x,y,0),$ for every $(x,y,t)\in \H$ (see Lemma 2.10); this inequality provides   information about the unitary ball w.r.t. the norm $N.$ Specifically, the following remark holds:
\begin{remark}\label{proiezione bolla norma}
Let $N:\H\to [0,+\infty)$ be a homogeneous norm. Moreover, let us suppose that $N$ is a radial H-convex function. Let $B_N(e,1)=\{\xi\in\H:\ N(\xi)\le 1\}$.
Then,
\begin{enumerate}[i.]
\item   $B_N(e,1)$  generates an H-convex family via dilation and the $\H$-cone-function with vertex $e$ and base $\partial B_N(e,1)$ is exactly $N;$
\item the previous mentioned lemma in \cite{BalFasSob2018} implies that
$$
{\rm Pr}_{\{t=0\}}(B_N(e,1))\subset \left(B_N(e,1)\cap \{(x,y,0)\in\H\}\right),
$$
where ${\rm Pr}_{\{t=0\}}$ denotes the projection onto the $xy$-plane.
\end{enumerate}
\end{remark}
As a matter of fact, there exist H-convex, homogeneous and radial functions that are not homogeneous norms (see, for instance, Example \ref{example non norma}).
In the following, we try to shed some light on the properties that a radial set $K$ satisfying \textbf{a.-c.} should enjoy in order to have some chance to generate an $\H$-cone-function.
%
%
%The question is if, given a  general radial set $K$ which generates an H-convex family via dilation,
%
%
%similar results can be inferred for a general radial set $K$ which generates an H-convex family via dilation. In this line of investigation, we are able to prove the following fundamental result concerning an \lq\lq interior constraint\rq\rq\ for the boundary of this set $K$:
\begin{lemma}\label{teo calotta sferica}
Let $K\subset \H,$ and let $\Kca$ be as in Lemma \ref{lemma solido di rivoluzione}. If $\xi_0=(x_0,y_0,t_0)\in K,$ $t_0>0,$ and $r_0=\|(x_0,y_0)\|_E>0$,
then $K$ contains the spherical cap of the Kor\'anyi ball $B^+_G,$ that is given by:
$$
B^+_G=\left\{\xi=(x,y,t)\in \H:\ \xi\in B_G\left(e,\|\xi_0\|_G\right),\ t\ge t_0\right\}.
$$
If $t_0<0$, then, symmetrically, $B^-_G\subset K,$ where $B^-_G$ is defined by setting $t\le t_0$ in the previous line.
\end{lemma}

\noindent\textbf{Proof:}  By the radiality of $K$, we restrict our attention to the profile of $K,$ i.e., the set $\partial K\cap \{y=0\}.$
If  $\xi_0=(x_0,y_0,t_0)\in K,$ then $\xi_0'=(r_0,0,t_0)\in K,$ where $r_0=\sqrt{x_0^2+y_0^2},$ and Lemma \ref{lemma solido di rivoluzione} implies that $S_{\xi_0}\subset K$. In particular, $K$ contains the image of the curve
$$
r\mapsto \Phi(r)=(r,0,\varphi(r)),\qquad r\in[r_0/\sqrt{2},r_0],$$
where $\varphi$ is given in Lemma \ref{lemma solido di rivoluzione}.
Since
$$
\varphi'(r_0)=-\frac{2r_0^3}{t_0}=-\frac{2r_0^3}{\varphi(r_0)},
$$
and since this argument can be repeated for every point in $K$ and, in particular, for all the points of the previous curve $\Phi$ inside $K$,
then, starting from $\xi_0'\in K,$ there exists a new curve $r\mapsto \Psi(r)=(r,0,\psi(r))\in K,$ with $\psi:[0,r_0]\to \R,$ that is a solution of the Cauchy problem
$$\left\{
\begin{array}{l}
\displaystyle \psi'(r)=-\frac{2r^3}{\psi(r)},\qquad r\in[0,r_0]\\
\displaystyle \psi(r_0)=t_0\\
\end{array}
\right.
$$
Note that this argument cannot be applied if $r >r_0$. The solution of the Cauchy problem above is, for $r\in[0,r_0]$,
$$
\psi(r)=\sqrt{-r^4+r^4_0+t_0^2},
$$
i.e. $(r^4+\psi^2(r))^{1/4}=(r^4_0+t_0^2)^{1/4}.$ Since the curve $\Psi$ belongs to $K$, and $K$ is H-convex and radial, we have the claim.

\QED

We are now ready to prove the first part of Theorem \ref{ball}:%give one of the main results of the paper:
%\begin{theorem}\label{ball}
% Let $K\subset \H$ be a radial and compact set which generates an H-convex family via dilations. If $\xi_0=(x_0,y_0,0)\in \partial K$, then
% $$
% B_G(e,\|\xi_0\|_G)\subset K.
% $$
%\end{theorem}
\psn\textbf{Proof of Theorem \ref{ball}-i}:
Without loss of generality, let us assume that $\xi_0=(x_0,0,0)\in \partial K$, with $x_0>0$.
By Proposition \ref{punti non caratteristici} the point $\xi_0$ is non-characteristic and the horizontal line $r$ through $\xi_0$ in \eqref{retta r} enters $K$. Clearly, there exists $\{\xi_n^+\}_{n> 0}\subset K$, with $\xi_n^+=(x_n,0,t_n), \ t_n> 0,\ \xi_n^+\to \xi_0$, and
there exists $\{\xi_n^-\}_{n> 0}\subset K$, with $\xi_n^-=(x_n,0,t_n), \ t_n< 0,\ \xi_n^-\to \xi_0$.
Finally, from Lemma \ref{teo calotta sferica}, taking the limit as $n\to \infty,$ by means of a continuity argument and the closedness of $K,$ we prove the claim.

\QED

The previous result has an immediate and nice application to homogeneous norms:
\begin{remark}
Let $N$ be a homogeneous and radial norm. Since $B_N(e,1)$ generates an H-convex family via dilation, then
$$
B_G(e,1)\subset B_N(e,1).$$
\end{remark}

Now we are ready to prove the second part of Theorem \ref{ball}, that generalizes Remark \ref{proiezione bolla norma}:
%
%\begin{theorem}\label{teorema proiezione}
% Let $K\subset \H$ be a radial and compact set which generates an H-convex family via dilations.
% Then, the projection of $K$ onto the plane  $\{t=0\}$ is contained in the intersection of $K$ with $\{t=0\}$, i.e.,
% $$
%{\rm Pr}_{\{t=0\}}(K)\subset \left(K\cap \{(x,y,0)\in\H\}\right).
%$$
%\end{theorem}
\psn\textbf{Proof of Theorem \ref{ball}-ii}: We argue by contradiction. Let $\overline{\xi}=(\overline{x},0,0)\in\partial K,$ with $\overline{x}>0,$ and suppose that $\xi_0=(x_0\overline{x},0,t_0)\in \partial K$ for some $x_0>1$ and $t_0>0$.
Since $\delta_{{1}/{\overline{x}}}\overline{\xi}=(1,0,0)\in \partial (\delta_{{1}/{\overline{x}}}K)$ and $\delta_{{1}/{\overline{x}}}\xi_0=(x_0,0,{t_0}/{\overline{x}^2})\in \partial (\delta_{1/{\overline{x}}}K)$, then, by
Remark \ref{remark dilatazioni}-b., there is no loss of generality if we take $\overline{x}=1$.

We will prove that the condition $\overline{\xi}=(1,0,0)\in\partial K$ and $\xi_0=(x_0,0,t_0)\in \partial K,$ for some $x_0>1$ and $t_0>0,$ leads to a contradiction.
From Theorem \ref{ball} we get that
$$
\xi_{t,\theta}=(r\cos \theta,r\sin\theta, t)\in \partial B_G(e,1)\subset K.
$$
with $r=\sqrt[4]{1-t^2},$ and for some fixed $-1<t<0$ and $\theta$.
Taking in \eqref{parabola 1 zero} $\xi_0=(x_0,0,t_0)\in K$,  the point $\xi_{t,\theta}$ can be obtained as the second point in \eqref{parabola 1 zero} by  considering
$$
\alpha=r\t\cos\theta-x_0,\qquad \beta=r\t\sin\theta,
$$
and $\t>0$ is a solution of the equation
\begin{equation}\label{eq:tau}
-t\t^2-2x_0r\t\sin\theta+t_0=0
\end{equation}
(we recall that $t<0$).
In order to do that, let us choose $\theta\in (0,\pi)$ so that $\sin\theta=\frac{\sqrt{-t_0 t}}{x_0r};$ for this value of $\theta$ we obtain $\t=\sqrt{\frac{t_0}{-t}}.$
In this case, the curve $\widetilde\xi$  in  Proposition \ref{prop condizione parabola}-i.  is given by
\begin{eqnarray*}
&&\widetilde\xi (s)= \Biggl(s,\frac{r(x_0-s)\sin\theta}{x_0-r\cos\theta},\\
&&\qquad\qquad \frac{t_0((\t-1)s-r\t\cos\theta +x_0)^2-2x_0r\t\sin\theta (x_0-s)((\t-1)s-r\t\cos\theta +x_0)}{\t^2(x_0-r\cos\theta)^2}
\Biggr),\qquad
\end{eqnarray*}
where $s\in [r\cos\theta,x_0].$
A tedious calculation gives that the curve $s\mapsto\widetilde\xi (s)$ intersects the $xy$-plane if and only if
$s=\overline{s}=\frac{t_0 r\cos\theta-\t tx_0}{t_0-\t t}.$ Let us prove the third component of $\xi_{t,\theta}$ can be chosen so that
$\|\widetilde\xi (\overline{s})\|_G>1$.
In fact, taking into account our choice of $r$ and $\theta$ (and hence of $\t$ defined by \eqref{eq:tau}),
\begin{eqnarray*}
\|\widetilde\xi (\overline{s})\|_G^2&=&\left(\frac{t_0 r\cos\theta-\t tx_0}{t_0-\t t}\right)^2+\left(\frac{t_0 r\sin\theta}{t_0 -\t t}\right)^2\\
&=&\frac{t_0 r^2+2\sqrt{t_0}\sqrt{-t}\sqrt{x_0^2r^2+t_0 t}-tx_0^2}{(\sqrt{t_0}+\sqrt{-t})^2}\\
&=&\frac{t_0 \sqrt{1-t^2}+2\sqrt{t_0}\sqrt{-t}\sqrt{x_0^2\sqrt{1-t^2}+t_0 t}-tx_0^2}{(\sqrt{t_0}+\sqrt{-t})^2}.
\end{eqnarray*}
If $t\to 0^-$ we have
$$
\|\widetilde\xi (s_0)\|_G^2=\frac{t_0+2x_0\sqrt{t_0}\sqrt{-t}+o(\sqrt{-t})}{t_0+2\sqrt{t_0}\sqrt{-t}+o(\sqrt{-t})}\to 1^+,
$$
since $t_0>0$ and $x_0>1$ are fixed.
Hence, if $t$ is negative and small enough, we find a point of the parabola  $\widetilde\xi$ that belongs to the $xy$-plane, but is not in $K$. This contradicts condition  \eqref{eq:condition CH} and the proof is finished.

\QED

The next example highlights that there exist H-convex radial homogeneous functions that are not homogeneous norms.
\begin{example}\label{example non norma}\rm{
Let us consider the function $f:\H\to \R$ defined as follows:
$$
f(x,y,t)=\begin{cases}
f_1(x,y,t),&\qquad t\ge 0\\
f_2(x,y,t),&\qquad t<0,
\end{cases}
$$
where $f_1,f_2:\H\to \R$ are the H-convex radial homogeneous functions given by
$$
f_1(x,y,t)=
\max\left\{\|(x,y)\|_E,\frac{1}{2}\|(x,y,t)\|_G\right\},\qquad f_2(x,y,t)=\|(x,y,t)\|_G.
$$
The function $f$ is trivially H-convex on the half-space $\overline{\H^+}=\{(x,y,t):\; t\ge 0\}$, as well as $\H^-=\{(x,y,t):\; t<0\}.$
Let us prove that $f$ is H-convex on the whole $\H,$ i.e., $f|_r$ is a convex function for every $r$, where $r$ denotes a horizontal line.
Note that, if $r\subset \overline{\H^+},$ or $r\subset \H^-,$ then there is nothing to prove. Thus, we can consider the horizontal lines $r$ having the parametric form
$$
s\mapsto (x_0+vs,y_0+v's,ws),\quad s\in \R,
$$
where $x_0,y_0, v,v',w\in \R,$ $(x_0,y_0)\neq (0,0),$ $(v,v')\neq (0,0),$ $w\neq 0.$  Since the function $f$ is radial, without loss of generality we can assume that $y_0=0,$ and set $x_0>0.$ Let us consider the restriction $f|_r:\R\to \R;$ if $w>0,$ then
$$
f|_r(s)=f(x_0+vs,v's,ws)=\begin{cases}
\max\left\{\|(x_0+vs, v's)\|_{E}, \frac{1}{2}\|(x_0+vs,v's,ws)\|_G\right\},&\qquad s\ge 0\\
\|(x_0+vs,v's,ws)\|_G,&\qquad s<0.
\end{cases}
$$
It is a simple exercise to prove that, if $g_1, g_2:\R\to \R$ are convex functions, and $g:\R\to \R$ is the function given by
$$
g(s)=\begin{cases}
g_1(s),&\quad s<0\\
g_2(s),&\quad s\ge 0,
\end{cases}
$$
then, under the assumptions that $g$ is continuous at $s=0,$ $D_-g(0),$ $D_+g(0)$ exist and $D_g(0)\le D_+g(0),$  the function  $g$ turns out to be convex on the whole $\R;$ this is true, in particular, if $Dg(0)$ is well defined.

In our case, the function $f$ is H-convex on $\overline{\H^+}$ and $\H^-,$ therefore $f|_r$ is convex on $(-\infty,0),$ and $[0,+\infty).$
Furthermore, $f|_r$ is differentiable at $s=0.$ Indeed, simple computations show that, whenever $s\to 0,$ then, if $v>0,$ $f|_r(s)=x_0+sv+o(s),$ while, if $v=0,$ $f|_r(s)=x_0+o(s).$ In both cases, $Df|_r(0)$ does exist. This proves that $f|_r$ is convex on $\R.$
The case $w<0$ can be treated in the same way.}
\end{example}

The following example is crucial: it shows that there exists a compact and radial set, H-convex, but not $\R^3$-convex, whose boundary is not the base for an $\H$-cone-function with vertex $e.$
\begin{example}\label{example importante}\rm{
Let us consider the function $f_1:\H\to \R$ defined by
$$f_1(x,y,t)=\left((x^2+y^2)^2+2+\frac{1}{2}\sin t\right)^{1/4}-2^{1/4}.$$
This is an H-convex function studied in \cite{BaRi2003}. Now let us define the constants $c$ and $c_0$ by
$$0<c_0=\sqrt[4]{(1+2^{1/4})^4-2}<c= \sqrt[4]{(1+2^{1/4})^4-3/2},$$
and consider the point $\overline \xi=(c,0,3\pi/2)$.
Define $f:\H\to \R$  by
$$
f(\xi)=\max\left\{f_1(\xi),\ \frac{\|\xi\|_G}{\|\overline\xi\|_G}\right\},
$$
and consider the set
$$
K=\lev_{\le 1}(f).
$$
The function $f$ is trivially H-convex, $f(e)=0$ and $f\bigl|_{\partial K}=1.$
Clearly, $e\in{\rm int}(K);$ moreover, $K$ is H-convex and compact. It is easy to see that
$$
\xi_0=(c_0,0,0)\in \partial H,\qquad \xi_1=(c_0,0,\pi/2)\not\in K,\qquad \xi_2=(c_0,0,\pi)\in\partial K;$$
therefore $K$ is not $\R^3$-convex. Finally,  we note that $\overline\xi\in \partial K$ and
$${\rm Pr}_{\{t=0\}}(\overline\xi)=(c,0,0)\not\in \left(K\cap \{(x,y,0)\in\H\}\right)=\{(x,y,0)\in\H:\ \|(x,y)\|_E\le c_0\};$$
hence, by Theorem \ref{ball},
  $K$ does not generate an H-convex family via dilation.
The previous arguments imply that there is neither an $\H$-cone-function, nor a cone-function, with vertex $e$ and base $\partial K.$
}
\end{example}
%
%
%
%\bibliography{capi}
%\bibliographystyle{plain}

\end{document}